\documentstyle[12pt,amscd,amssymb]{amsart}

\def\sw#1{{\sb{(#1)}}}
\def\tens{\mathop{\otimes}}

\def\<{{\langle}}
\def\>{{\rangle}}

\def\eps{\epsilon}

\def\note#1{{}}
\def\can{{\rm can}}
\def\note#1{}
\def\M{{\bf M}}
\def\z{{\frak z}}
\def\y{{\frak y}}

\def\tA{{\tilde{A}}}
\def\tC{{\tilde{C}}}
\def\tD{{\tilde{D}}}
\def\tpsi{{\tilde{\psi}}}
\def\tmu{{\mu_\tA}}
\def\tnu{{\tilde{\nu}}}
\def\teps{{\epsilon_\tC}}
\def\tDelta{{\Delta_\tC}}
\def\tM{{\tilde{M}}}
\def\ta{{\tilde{a}}}
\def\tc{{\tilde{c}}}
\def\tm{{\tilde{m}}}
\def\coker{{\rm coker}}
\def\cocan{{\rm cocan}}
\def\eq{{\rm eq}}

\def\twist{{\rm twist}}

\def\Label{\label}
\begin{document}

\newtheorem{proposition}{Proposition}[section]
\newtheorem{lemma}[proposition]{Lemma}
\newtheorem{corollary}[proposition]{Corollary}
\newtheorem{theorem}[proposition]{Theorem}

\theoremstyle{definition}
\newtheorem{definition}[proposition]{Definition}
\newtheorem{example}[proposition]{Example}

\theoremstyle{remark}
\newtheorem{remark}[proposition]{Remark}

\newcommand{\Section}{\setcounter{definition}{0}\section}
\newcounter{c}
\renewcommand{\[}{\setcounter{c}{1}$$}
\newcommand{\etyk}[1]{\vspace{-7.4mm}$$\begin{equation}\Label{#1}
\addtocounter{c}{1}}
\renewcommand{\]}{\ifnum \value{c}=1 $$\else \end{equation}\fi}

\title[Coalgebra-Galois extensions]{Coalgebra-Galois extensions from the
extension theory point of view}
\author{Tomasz Brzezi\'nski}
\address{Department of Mathematics, University of York, Heslington, 
York YO10 5DD, U.K.}
\email{tb10@@york.ac.uk}
\urladdr{http//www.york.ac.uk/\~{}tb10}
\thanks{Lloyd's of London Tercentenary Fellow.}
\thanks{On leave from: Department of Theoretical Physics, University of
\L\'od\'z, Pomorska 149/153, 90-236, \L\'od\'z, Poland}
\begin{abstract}
Coalgebra-Galois extensions generalise Hopf-Galois extensions, which can
be viewed as non-commutative torsors. In this
paper it is analysed when a coalgebra-Galois extension is a separable,
split, or strongly separable extension.
\end{abstract}
\maketitle
\section{Introduction}
\noindent Given a coalgebra $C$, an algebra $A$ and a right coaction
$\rho^A: A\to A\otimes C$ one can define a fixed point subalgebra $B$ of
$A$ as consisting of all those elements of $A$ over which the coaction is
left-linear. In this way one obtains an extension $B\hookrightarrow A$,
which is called a {\em coalgebra-Galois extension} if a certain canonical
left $A$-module, right  $C$-comodule map is bijective \cite{BrzMa:coa}
\cite{BrzHaj:coa}.  The aim of this article is to analyse such
coalgebra-Galois extensions from the extension
theory point of view. In particular we study the problem when such
extensions are separable, split or strongly separable extensions. This
problem is put in a broader context of {\em entwining structures} and
{\em entwined modules} introduced in \cite{BrzMa:coa} \cite{Brz:mod},
as a generalisation of a Doi-Hopf datum and Doi-Koppinen
modules \cite{Doi:uni} \cite{Kop:var}, respectively. We make  use of
the notion of a separability of a functor introduced in \cite{Nas:sep},
and, as a byproduct, we generalise some of the results of
\cite{CaeMil:sep} obtained recently for Doi-Koppinen modules. 

The paper is organised as follows. In Section~2 we recall definitions
 and give examples of entwining structures and 
entwined modules. In Section~3 we analyse when certain functors between
categories of entwined modules induced 
by morphisms of entwining structures are separable. In Section~4 we
 apply the
results of Section~3 to prove  that a sufficient and
necessary condition for a coalgebra-Galois
extension to be  separable is 
the separability of a certain induction functor. This, in turn, is
equivalent to  the existence of a
normalised integral in the canonical entwining structure. In Section~5
we analyse when a coalgebra-Galois extension is a split extension. This
turns out to be related to the separability of the forgetful functor
from the category of entwined modules to the category of right modules -
another special case of the main theorem in Section~3.
Finally, in Section~6 we study the problem when a coalgebra-Galois 
extension is a
strongly separable extension in the sense of \cite{Kad:Jon}.

We work over a commutative ring $k$ with identity 1. We  assume that all
the algebras are over $k$, associative  and  unital, and the coalgebras 
are over $k$, coassociative 
and  counital. Unadorned tensor product is over $k$. For any 
$k$-modules 
$V, W$ the symbol  ${\rm Hom}(V,W)$ denotes the $k$-module 
of $k$-linear maps $V\to W$ and the identity map $V\to V$ is
 denoted by $V$. The twist map between $k$-modules $V, W$ is denoted by
$\twist: V\otimes W\to W\otimes V$, $v\otimes
w\mapsto w\otimes v$. We also implicitly identify $V$ 
with $V\otimes k$
and $k\otimes V$ via the canonical isomorphisms. 

For a $k$-algebra $A$ we use $\mu_A$ to denote the product as a map and 
$1_A$ to denote  the identity both as an element of $A$ and as a map $k\to
A$, $\alpha\mapsto \alpha 1_A$. $\M_A$ (resp.\ ${}_A\M$) denotes 
the category of right (resp.\ left) $A$-modules. The morphisms in this
category are denoted by ${\rm Hom}_A(M,N)$ (resp.\ ${}_A{\rm Hom}(M,N)$).
For any $M\in \M_A$ (resp.\ $M\in {}_A\M$), the symbol $\rho_M$ (resp.\
${}_M\rho$) denotes the action as a map (on elements the action  is
denoted by a dot). We often write $M_A$ (resp.\ ${}_AM$) to
indicate in which context the $A$-module $M$ appears.
For any $M\in \M^A$, $N\in {}^A\M$ we will write $\eq_{M_AN}:M\otimes
A\otimes N\to M\otimes N$ for the action equalising map defining 
tensor product $M\otimes_A N$, i.e., $\eq_{M_AN} = \rho_M\otimes
N-M\otimes{}_N\rho$, $M\otimes_A N = \coker(\eq_{M_AN})$.

For a
$k$-coalgebra $C$ we use $\Delta_C$ to denote the coproduct and $\eps_C$ to
denote the counit. Notation for comodules is similar to that for modules
but with subscripts replaced by superscripts, i.e. $\M^C$ is the
category of right $C$-comodules, $\rho^M$ is a right coaction etc. We
use the Sweedler notation for coproducts and coactions, i.e. $\Delta_C(c)
= c\sw 1\otimes c\sw 2$, $\rho^M(m) = m\sw 0\otimes m\sw 1$ (summation
understood).   For any $V \in \M^C$, 
$W\in {}^C\M$,     $V \square_{C}W$ denotes the cotensor product, which 
is defined by the exact sequence    
\[    
\begin{CD}
0 @>>> V\square_{C}W @>>> V\otimes W
@>{\eq^{V^CW}}>> V \otimes C\otimes W,  
\end{CD} 
\]    
where $\eq^{V^CW}$ is the coaction equalising map, i.e., $\eq^{V^CW} = \rho^{V}\otimes W -
  V\otimes {}^{W}\!\rho$.

\section{Preliminaries on entwining structures and coalgebra-Galois
extensions} 
\begin{definition}
\Label{ent}
An {\em  entwining
structure} (over $k$) is a 
triple $(A,C)_\psi$ consisting of a $k$-algebra $A$, a $k$-coalgebra
$C$ and a $k$-linear map $\psi: C\tens A\to A\tens C$ satisfying
\begin{equation}
\psi\circ(C\tens \mu_A) = (\mu_A\tens C)\circ (A\tens\psi)\circ(\psi\tens A),
\quad \psi\circ (C\tens 1_A) = 1_A\tens C,
\label{diag.A}
\end{equation}
\begin{equation}
(A\tens\Delta_C)\circ\psi = (\psi\tens
C)\circ(C\tens\psi)\circ(\Delta_C\tens A), \quad (A\tens \eps_C)\circ\psi =
\eps_C\tens A.
\label{diag.B}
\end{equation}
A morphism of entwining structures is a pair
 $(f,g):(A,C)_\psi\to
(\tilde{A}, \tilde{C})_{\tilde{\psi}}$, where $f:A\to \tilde{A}$ is
an algebra map, $g: C\to \tilde{C}$ is a coalgebra map, and
$(f\otimes g)\circ\psi =
\tilde{\psi}\circ(g\otimes f)$. 
\end{definition}
The category of entwining structures is a tensor category
with tensor product $(A,C)_\psi \otimes (\tA, \tC)_\tpsi =
(A\otimes\tA, C\otimes\tC)_{(A\otimes\twist\otimes \tC)\circ
(\psi\otimes\tpsi)\circ (C\otimes\twist\otimes \tA)}$, and unit object
$(k,k)_\twist$. 

For $(A,C)_\psi$ we  use the
notation $\psi(c\otimes a) = a_\alpha \otimes c^\alpha$ (summation
over a Greek index understood), for all $a\in
A$, $c\in C$.
The notion of an entwining structure was  introduced in
\cite[Definition~2.1]{BrzMa:coa}. It is self-dual in the sense that
conditions in Definition~\ref{ent} are invariant under the operation
consisting of  interchanging
of $A$ with $C$, $\mu_A$ with $\Delta_C$, and  $1_A$ with $\eps_C$, and
reversing the order of maps. Below are two classes of examples
of entwining  structures coming from Galois-extensions. 

\begin{example}[\cite{BrzHaj:coa}]
Let $C$ be a coalgebra, $A$ an algebra and a right $C$-comodule. Let 
$B:= \{b\in A\; | \; \rho^A(ba) =
b\rho^A(a)\}$ and assume that the canonical left $A$-module, right
$C$-comodule map
$
\can:A\otimes _BA\to A\otimes C$, $a\otimes a'\mapsto a\rho^A(a')$,
is bijective. Let $\psi:C\otimes A\to A\otimes C$ be a $k$-linear map
given by
$
\psi(c\otimes a) = \can(\can^{-1}(1_A\otimes c)a).
$
Then $(A,C)_\psi $ is an entwining structure. The 
extension $B\hookrightarrow A$ is called a
{\em coalgebra-Galois extension} (or a {\em $C$-Galois extension})  and is
denoted by $A(B)^C$. $(A,C)_\psi$ is the
{\em canonical entwining structure} associated to $A(B)^C$. A
coalgebra-Galois extension $A(B)^C$ is said to be {\em copointed} if
there exists a group-like $e\in C$ such that $\rho^A(1_A)=1_A\otimes e$.
\label{can.ex}
\end{example} 
Dually we have
\begin{example}[\cite{BrzHaj:coa}]
Let $A$ be an algebra, $C$ a coalgebra and a right $A$-module. Let 
$B:= C/I$, where $I$ is a coideal in $C$, 
$$
I:={\rm span}\{(c\cdot a)\sw 1\xi((c\cdot a)\sw 2)-c\sw 1\xi(c\sw
2\cdot a)\; |\; a\in A,\; c\in C,\; \xi\in C^*\},
$$
 and assume that the canonical left $C$-comodule, right
$A$-module map
$
\cocan:C\otimes A\to C\square_B C$, $c\otimes a \mapsto c\sw 1\otimes c\sw
2\cdot a$,
is bijective. Let $\psi:C\otimes A\to A\otimes C$ be a $k$-linear map
given by
$
\psi = (\eps_C\otimes A\otimes C)\circ (\cocan^{-1}\otimes C)\circ
(C\otimes \Delta_C)\circ \cocan.
$
Then $(A,C)_\psi $ is an entwining structure. The 
coextension $C\twoheadrightarrow B$ is called an
{\em algebra-Galois coextension} (or an {\em $A$-Galois coextension})  
and is
denoted by $C(B)_A$. $(A,C)_\psi$ is the
{\em canonical entwining structure} associated to $C(B)_A$. An
algebra-Galois coextension $C(B)_A$ is said to be {\em pointed} if there
exists an algebra map $\kappa:A\to k$ such that $\eps_C\circ \rho_C =
\eps_C\otimes\kappa$. 
\label{cocan.ex}
\end{example} 

Associated to an entwining structure is the category of entwined
modules. 
\begin{definition}
Let $(A,C)_\psi$ be an entwining structure. An (entwined) 
{\em $(A,C)_\psi$-module}
is a right $A$-module, right $C$-comodule $M$ 
such that 
$$
\rho^M\circ\rho_M = (\rho_M\tens C)\circ(M\tens\psi)\circ(\rho ^M\otimes
A), 
$$
(explicitly: 
$
\rho^M(m\cdot a) = m\sw 0\cdot a_\alpha\tens m\sw 1^\alpha$, $ \forall
a\in A, m\in M$).
A morphism of $(A,C)_\psi$-modules is a right $A$-module map
which is also a right $C$-comodule map. The category of
$(A,C)_\psi$-modules is denoted by $\M_A^C(\psi)$.
\label{def.mpsi}
\end{definition}
The category $\M_A^C(\psi)$ was introduced and
studied in \cite{Brz:mod}. An example of such modules are Doi-Koppinen
modules introduced in \cite{Doi:uni}, \cite{Kop:var}. 
In this paper we will be concerned with two covariant functors between
categories of entwined modules, which are
special cases of the construction in 
\cite[Section~3]{Brz:mod}\footnote{Although the paper \cite{Brz:mod} is
restricted to $k$ being a field, all the results of \cite{Brz:mod}
quoted in the present paper can easily be seen to hold for a general
$k$.} (see also \cite{CaeRai:ind} for the Doi-Koppinen case). 
These functors are induced by certain
morphisms of entwining structures.
\begin{definition}
Let $(f,g) :(A,C)_\psi \to (\tA,\tC)_\tpsi$ be a morphism of 
entwining structures. View  $C$ as a left $\tC$-comodule via
${}^C\!\rho = (g\otimes C)\circ\Delta_C$ and $C\otimes\tA$ as a right
$\tC$-comodule via
$\rho^{C\otimes\tA} = (C\otimes\tpsi)\circ(C\otimes
g\otimes\tA)\circ(\Delta_C\otimes\tA)$. Then $(f,g)$ is said to be an
{\em admissible
morphism} iff:\\
\indent (i) for all $\tM\in\M_\tA^\tC(\tpsi)$,
$\tM\square_\tC (C\otimes C) = (\tM\square_\tC C)\otimes C$,\\
\indent (ii) for all $M\in \M_A$, $(M\otimes
(C\otimes\tA))\square_\tC C = M\otimes
((C\otimes\tA)\square_\tC C)$.
\label{admissible}
\end{definition}
For example, if $C,\tC$ are $k$-flat then $(f,g)$ is
an admissible morphism provided that ${}^\tC \!C$ is coflat. On the
other hand if $k$ is a regular 
ring or a field every
morphism is admissible. Also, it can be easily checked that
 the following morphisms 
$(A,\eps_C): (A,C)_\psi\to (A,k)_\twist$ and $(1_A,
C): (k,C)_\twist\to (A,C)_\psi$ are admissible.  
\begin{example}
Let $(f,g) :(A,C)_\psi \to (\tA,\tC)_\tpsi$ be an admissible morphism of 
entwining structures.
View $\tA$ as a right $A$-module via $\rho_\tA = \tmu\circ(\tA\otimes
f)$, and $C$ as a right $\tC$-comodule via
$\rho^C = (C\otimes g)\circ\Delta_C$. Then:

(1) For any $\tM\in\M_\tA^\tC(\tpsi)$, $\tM\square_\tC C$ is an
$(A,C)_\psi$-module with structure maps $\rho^{\tM\square_\tC C} =
\tM\otimes \Delta_C$ and
$$
\rho_{\tM\square_\tC C} :\tM\square_\tC C\otimes A\to \tM\square_\tC C, \qquad
\sum_i \tm_i\otimes  c_i\otimes a =
\sum_i\tm_if(a_\alpha)\otimes c_i^\alpha.
$$

(2) For any $M\in\M_A^C(\psi)$, $M\otimes_A\tA$ is an
$(\tA,\tC)_\tpsi$-module with structure maps $\rho_{M\otimes_A\tA} =
M\otimes_A\tmu$
and
$$
\rho^{M\otimes_A\tA}:M\otimes_A\tA\to M\otimes_A\tA\otimes\tC, \qquad
m\otimes\ta\mapsto m\sw0\otimes\ta_\alpha\otimes g(m\sw 1)^\alpha.
$$

(3) The covariant functor $-\square_\tC C:\M_\tA^\tC(\tpsi) \to
\M_A^C(\psi)$ is the right adjoint of $-\otimes_A\tA:\M_A^C(\psi)\to
\M_\tA^\tC(\tpsi)$. The adjunctions are:
$$
\forall M\in\M_A^C(\psi), \qquad \Phi_M: M\to (M\otimes_A\tA)\square_\tC C,
\qquad m\mapsto m\sw 0\otimes 1_\tA\otimes m\sw 1,
$$
$$
\forall \tM\in \M_\tA^\tC(\tpsi), \quad \Psi_\tM : (\tM\square_\tC
C)\otimes _A\tA\to \tM, \quad \sum_i \tm_i\otimes c_i\otimes \ta
\mapsto \sum_i\tm_i\cdot\ta\eps_C(c_i).
$$
\label{ex.functors}
\end{example}
Applying
Example~\ref{ex.functors} to morphisms  
$(A,\eps_C): (A,C)_\psi\to (A,k)_\twist$ and $(1_A,
C): (k,C)_\twist\to (A,C)_\psi$ one obtains

\begin{example} Let $(A,C)_\psi$ be an entwining structure. Then 

(1) If $M$ is a right $A$-module then $M\otimes C$ is an
$(A,C)_\psi$-module with
the coaction $M\otimes \Delta_C$ and the action $(m\otimes c)\cdot a =
m\cdot\psi(c\otimes a)$, for all $a\in A, c\in C$ and $m\in M$. In
particular $A\otimes C\in \M^C_A(\psi)$. The operation
$M\mapsto M\otimes C$ defines a 
covariant functor $-\otimes C:\M_A\to \M_A^C(\psi)$ which 
is the right adjoint of the forgetful functor $\M_A^C(\psi)\to \M_A$.
\note{The adjunctions are
$$
\forall M\in \M_A^C(\psi), \quad \Phi_M: M\to M\otimes C, \quad \Phi_M =
\rho^M, 
$$
$$
\forall M\in \M_A, \quad \Psi_M:M\otimes C\to M, \quad \Psi_M =
M\otimes\eps_C. 
$$}

(2) If $V$ is a right $C$-comodule then $V\otimes A\in \M^C_A(\psi)$
with the action $V\otimes\mu_A$ and the coaction $v\otimes a\mapsto v\sw
0\otimes \psi(v\sw 1\otimes a)$ for any $a\in A$ and $v\in V$. In
particular $C\otimes A\in \M^C_A(\psi)$. The operation
$V\mapsto V\otimes A$ defines a 
covariant functor $-\otimes C:\M^C\to \M_A^C(\psi)$, which is the left 
adjoint of the forgetful functor $\M_A^C(\psi)\to \M^C$. 
\note{The adjunctions
are:
$$
\forall V\in \M^C, \quad \Phi_V: V\to V\otimes A, \quad \Phi_V =
V\otimes 1_A,
$$
$$\forall M\in \M_A^C(\psi), \quad \Psi_M: M\otimes A\to M, \quad \Psi_M =
\rho_M.
$$}
\label{ex.psi.modules}
\end{example}
Another class of examples of entwined modules comes from 
(co)algebra-Galois (co)extensions \cite{BrzHaj:coa}

\begin{example}
~\\
\indent (1) Let $(A,C)_\psi$ be the canonical entwining structure associated to
a coalgebra-Galois extension $A(B)^C$. Then $A$ is an
$(A,C)_\psi$-module via $\rho^A$ and $\mu_A$.

(2) Let $(A,C)_\psi$ be the canonical entwining structure associated to
an algebra-Galois coextension $C(B)_A$. Then $C$ is an
$(A,C)_\psi$-module via $\rho_C$ and $\Delta_C$.
\label{ex.psi.galois}
\end{example}

\section{Separable functors of entwined modules}
\noindent In this section we analyse when functors described in
Example~\ref{ex.functors} are separable. Recall from \cite{Nas:sep} that
a covariant functor $F :{\cal C}\to {\cal D}$ is {\em separable} if the
natural transformation ${\rm Hom}_{\cal C}(-,-)\to {\rm Hom}_{\cal
D}(F(-), F(-))$ splits. In this paper we are dealing with the pairs of
adjoint functors, so that the following characterisation of separable
functors, obtained in \cite{Raf:sep} \cite{Rio:cat}, is of  great 
importance
\begin{theorem}
Let $G:{\cal D}\to {\cal C}$ be the right adjoint of $F :{\cal C}\to {\cal
D}$ with adjunctions $\Phi: 1_{\cal C}\to GF$ and $\Psi:FG\to 1_{\cal
D}$. Then

(1) $F$ is separable if and only if $\Phi$ splits, i.e., for all objects
$C\in{\cal C}$ there exists a morphism $\nu_C \in{\rm Mor}_{\cal C}
(GF(C), C)$ such that $\nu_C\circ \Phi_C = C$ and for all $f\in {\rm
Mor}_{\cal C} (C, \tC)$, $\nu_\tC\circ GF(f) = f\circ\nu_C$.

(2) $G$ is separable if and only if $\Psi$ cosplits, i.e., for all
objects $D\in {\cal D}$ there exists a morphism $\nu_D \in{\rm Mor}_{\cal D}
(D, FG(D))$ such that $\Psi_D\circ \nu_D = D$ and for all $f\in {\rm
Mor}_{\cal D} (D, \tD)$, $\nu_\tD\circ f = FG(f)\circ\nu_D$. 
\label{thm.sep}
\end{theorem}

\begin{definition}
An admissible  morphism $(f,g): (A,C)_\psi\to (\tA,\tC)_\tpsi$ of entwining
structures is said to be:

(1) {\em integrable} if there exists $\lambda \in {\rm Hom}_A((C\otimes
\tA)\square_\tC C, A)$ such that the following diagrams commute
\begin{equation}
\begin{CD}
(C\otimes A\otimes \tA)\square_\tC C @>{\psi\otimes\tA\otimes C}>>
A\otimes (C\otimes \tA)\square_\tC C @>{A\otimes\lambda}>> A\otimes A\\
@VV{C\otimes f\otimes\tA\otimes C}V     @.     @VV{\mu_A}V\\
(C\otimes \tA\otimes \tA)\square_\tC C @>{C\otimes\tmu\otimes C}>>
(C\otimes \tA)\square_\tC C @>{\lambda}>> A,
\end{CD}
\label{int.a}
\end{equation}
\begin{equation}
\begin{CD}
(C\otimes\tA)\square_\tC C @>{C\otimes\tA\otimes \Delta_C}>> 
(C\otimes \tA)\square_\tC C\otimes C @>{\lambda\otimes C}>> A\otimes C\\
@VV{\Delta_C\otimes\tA\otimes C}V   @.   @|\\
C\otimes (C\otimes \tA)\square_\tC C @>{C\otimes\lambda}>> C\otimes A
@>{\psi}>> A\otimes C.
\end{CD}
\label{int.b}
\end{equation}
The right $A$-module structure of $(C\otimes\tA)\square_\tC C$ is as in
Example~\ref{ex.functors}(1), explicitly $\rho_{(C\otimes\tA)\square_\tC C}:
c'\otimes \ta\otimes c\otimes a\mapsto c'\otimes \ta f(a_\alpha)\otimes
c^\alpha$. 

(2) {\em totally integrable}, if there exists 
$\lambda\in {\rm Hom}_A((C\otimes
\tA)\square_\tC C, A)$ making it an integrable morphism and such that the
following diagram
\begin{equation}
\begin{CD}
C @>{\Delta_C}>> C\otimes C\\
@VV{1_A\circ\eps_C}V   @VV{C\otimes 1_\tA\otimes C}V\\
A @<{\lambda}<<  (C\otimes \tA)\square_\tC C
\end{CD}
\label{int.c}
\end{equation}
commutes.
\label{def.int.morphism}
\end{definition}
Notice that the condition (\ref{int.a}) makes sense because $\psi$ is a
morphism in $\M_A^C(\psi)$, $(f,g)$ is admissible and $(C\otimes\tmu)\circ(C\otimes f\otimes
\tA): C\otimes A\otimes\tA\to C\otimes\tA$ is a left $\tC$-comodule
map, where the $k$-modules involved are left $\tC$-comodules via
$(g\otimes C)\circ\Delta_C \otimes A\otimes\tA$ and $(g\otimes C)\circ\Delta_C \otimes
\tA$, respectively. Similarly, condition (\ref{int.b}) makes
sense because $\Delta_C$ is a left $\tC$-comodule map and
$\Delta_C\otimes\tA$ is a morphism in $\M_\tA^\tC(\tpsi)$. Dually to
Definition~\ref{def.int.morphism}
one considers
\begin{definition}
An admissible morphism $(f,g): (A,C)_\psi\to (\tA,\tC)_\tpsi$ of entwining
structures is said to be:

(1) {\em cointegrable} if there exists $\z \in {\rm Hom}^\tC(\tC, 
(\tA\otimes
C)\otimes_A \tA)$ such that the following diagrams commute
\begin{equation}
\begin{CD}
\tC @>{\z}>> (\tA\otimes C)\otimes_A\tA
@>{(\tA\otimes\Delta_C)\otimes_A\tA}>> (\tA\otimes C\otimes C)\otimes_A\tA
\\ 
@VV{\tDelta}V    @.     @VV{(\tA\otimes g\otimes C)\otimes_A\tA}V\\
\tC\otimes \tC @>{\tC\otimes \z}>> \tC\otimes (\tA\otimes C)\otimes_A\tA
@>{\tpsi\otimes\tC\otimes\tA}>> (\tA\otimes\tC\otimes C)\otimes_A\tA
\end{CD}
\label{coint.a}
\end{equation}
\begin{equation}
\begin{CD}
\tC\otimes\tA @>{\tpsi}>> \tA\otimes\tC @>{\tA\otimes \z}>> (\tA\otimes
\tA\otimes C)\otimes_A\tA\\
@|   @.   @VV{(\tmu\otimes C)\otimes_A\tA}V\\
\tC\otimes\tA @>{\z\otimes\tA}>> (\tA\otimes C)\otimes_A\tA\otimes\tA
@>{(\tA\otimes C)_A\otimes\tmu}>> (\tA\otimes C)\otimes_A\tA.
\end{CD}
\label{coint.b}
\end{equation}
The right $\tC$-comodule structure of $(\tA\otimes C)\otimes_A\tA$ is as
in
Example~\ref{ex.functors}(2), explicitly 
$\rho^{(\tA\otimes C)\otimes_A\tA}:
\ta\otimes c\otimes\ta' \mapsto \ta\otimes c\sw 1\otimes \ta'_\alpha
\otimes g(c\sw 2)^\alpha$. 

(2) {\em totally cointegrable}, if there exists $\z \in {\rm Hom}^\tC(\tC, 
(\tA\otimes C)\otimes_A \tA)$ making it a cointegrable morphism and such 
that the following diagram
\begin{equation}
\begin{CD}
\tC @>{\z}>> (\tA\otimes C)\otimes_A\tA\\
@VV{1_\tA\circ\teps}V   @VV{(\tA\otimes\eps_C)\otimes_A\tA}V\\
\tA @<{\mu_{\tA A}}<<  \tA\otimes_A\tA
\end{CD}
\label{coint.c}
\end{equation}
commutes. Here $\mu_{\tA, A}:\tA\otimes_A\tA\to \tA$ is the natural map induced
by $\mu_\tA$.
\label{def.coint.morphism}
\end{definition}

The right actions of $A$ on the $k$-modules involved in the above definition
are as follows. For any $a\in A$, $\ta,\ta'\in\tA$, $c,c'\in C$, $\tc\in
\tC$: $(\ta\otimes c)\cdot a = \ta f(a_\alpha)\otimes c^\alpha$,
$(\ta\otimes c\otimes c') = \ta f(a_{\alpha\beta})\otimes c^\beta\otimes
c^{\prime\alpha}$, $(\ta\otimes \tc\otimes c)\cdot a = \ta
f(a_\alpha)_\beta\otimes \tc^\beta\otimes c^\alpha$, $(\tc\otimes
\ta\otimes c)\cdot a = \tc\otimes \ta f(a_\alpha)\otimes c^\alpha$,
$(\ta\otimes\ta'\otimes c)\cdot a  = \ta\otimes \ta' f(a_\alpha)\otimes
c^\alpha$. Using properties of entwining structures and the fact that
$(f,g)$ is a morphism of entwining structures one can easily convince
oneself that all the maps featuring in
Definition~\ref{def.coint.morphism} are well-defined. 

With these definitions at hand we can now state the main result of this
section.
\begin{theorem}
Let $(f,g): (A,C)_\psi\to (\tA,\tC)_\tpsi$ be an admissible 
 morphism of entwining structures. 

(1) If for all $M\in \M_A^C(\psi)$, 
$(M\otimes_A\tA)\square_\tC C \subseteq
\coker(\eq_{M_A\tA}\square_\tC C)$,
then the functor 
$-\otimes_A\tA:\M_A^C(\psi)\to \M_\tA^\tC(\tpsi)$ is
separable if and only if $(f,g)$ is totally integrable.

(2) If for all $\tM\in\M_\tA^\tC(\tpsi)$,
$\ker(\eq^{\tM^\tC C}\otimes_A\tA)
\subseteq 
(\tM\square_\tC C)\otimes_A\tA$, then
 the functor $-\square_\tC C:\M_\tA^\tC(\tpsi)\to \M_A^C(\psi)$ is
separable if and only if $(f,g)$ is totally cointegrable.
\label{thm.main}
\end{theorem}
\begin{proof}
(1) Let $(f,g)$ be totally integrable and assume that $\lambda$ is as in
Definition~\ref{def.int.morphism}. For all $M\in \M_A^C(\psi)$ define
$$
\tnu_M : (M\otimes\tA)\square_\tC C\to M, \qquad
\sum_im_i\otimes\ta_i\otimes c_i \mapsto \sum_i m_i\sw
0\cdot\lambda(m_i\sw 1\otimes\ta_i\otimes c_i).
$$
Notice that the map $\tnu_M$ is well-defined since the fact that $(f,g)$
is admissible implies that for any
$x\in (M\otimes\tA)\square_\tC C$, one has $(\rho^M\otimes \tA\otimes
C)(x)\in (M\otimes (C\otimes\tA))\square_\tC C) =  M\otimes 
((C\otimes\tA)\square_\tC C)$. Take any 
$x= \sum_im_i\cdot
a_i\otimes \ta_i\otimes c_i\in (M\otimes\tA)\square_\tC C$. Then
\begin{eqnarray*}
\tnu_M(x) & = & \sum_i(m_i\cdot a_i)\sw 0\cdot\lambda((m_i\cdot a_i)\sw
1\otimes \ta_i\otimes c_i)\\
& = &  \sum_im_i\sw 0\cdot a_{i\alpha}\lambda(m_i\sw
1^\alpha \otimes \ta_i\otimes c_i)\qquad \qquad \mbox{\rm
($M\in\M_A^C(\psi)$)} \\
& = &  \sum_im_i\sw 0\cdot \lambda(m_i\sw
1 \otimes f(a_i)\ta_i\otimes c_i)\qquad \qquad \mbox{\rm
(by (\ref{int.a}))} \\
& = & \tnu_M(\sum_i m_i\otimes f(a_i)\ta_i\otimes c_i).
\end{eqnarray*}
The above calculation means that ${\rm Im}(\eq_{M_A\tA}\square_\tC
C)\subseteq \ker\tnu_M$, and together with the assumption that 
$-\square_\tC C$ preserves
the cokernel of the action equalising map 
$\eq_{M_A\tA}$ imply
that one can define the map $\nu_M:(M\otimes_A\tA)\square_\tC C\to
M$ by the diagram
\[
\begin{CD}
(M\otimes_A \tA)\square_\tC C \!@>>> \!\coker(\eq_{M_A\tA}\square_\tC
C)\!  @>>>
\!((M\otimes \tA)\square_\tC C)/\ker\tnu_M \! @>>>\!  M\\
@. @. @AAA @|\\
@. @. (M\otimes \tA)\square_\tC C @>{\tnu_M}>> M
\end{CD}
\]
 Slightly abusing the notation we will still write $\nu_M: 
\sum_im_i\otimes\ta_i\otimes c_i \mapsto \sum_i m_i\sw
0\cdot\lambda(m_i\sw 1\otimes\ta_i\otimes c_i)$.

To show that $\nu_M$ is a right $A$-module map, take any $a\in A$ and 
$x=\sum_i
m_i\otimes \ta_i\otimes c_i \in (M\otimes_A\tA)\square_\tC C$ and compute
\begin{eqnarray*}
\nu_M(x\cdot a) & = & \nu_M(\sum_im_i\otimes \ta_if(a_\alpha)\otimes
c_i^\alpha) \\
& = & \sum_i m_i\sw 0\cdot \lambda(m_i\sw 1\otimes
\ta_if(a_\alpha)\otimes c_i^\alpha)\\
& = & \sum_i m_i\sw 0\cdot \lambda(m_i\sw 1\otimes
\ta_i\otimes c_i)a \qquad \mbox{\rm ($f\in {\rm Hom}_A((C\otimes
\tA)\square_\tC C, A)$)}\\
& = & \nu_M(x)\cdot a.
\end{eqnarray*}
Furthermore we have
\begin{eqnarray*}
\nu_M(x\sw 0)\otimes x\sw 1 & = & \sum_i m_i\sw 0\cdot 
\lambda(m_i\sw 1\otimes
\ta_i\otimes c_i\sw 1)\otimes c_i\sw 2\\
& = &\sum_i m_i\sw 0\cdot \lambda(m_i\sw 2\otimes
\ta_i\otimes c_i)_\alpha\otimes m_i\sw 1^\alpha \qquad \mbox{\rm (by
(\ref{int.b}))} \\
& = &  \sum_i \rho^M(m_i\sw 0\cdot \lambda(m_i\sw 1\otimes
\ta_i\otimes c_i)) \qquad \qquad \mbox{\rm ($M\in \M_A^C(\psi)$)} \\
& = &  \rho^M\circ\nu_M(x),
\end{eqnarray*}
which proves that $\nu_M$ is a right $C$-comodule map. Using
(\ref{int.c}) one easily finds that the adjunction $\Phi_M$ is splitted by $\nu_M$. It
remains to be shown that $\nu_M$ is natural in $\M_A^C(\psi)$. Take any $M,N\in
\M_A^C(\psi)$ and $\phi \in {\rm Hom}_A^C(M,N)$. Then
\begin{eqnarray*}
\nu_N(\sum_i\phi(m_i)\otimes \ta_i\otimes c_i) &=& \sum_i \phi(m_i)\sw
0\cdot\lambda(\phi(m_i)\sw 1\otimes\ta_i\otimes c_i)\\
&=& \sum_i \phi(m_i\sw 0)\cdot\lambda(m_i\sw 1\otimes\ta_i\otimes c_i)
\\
&=& \sum_i \phi(m_i\sw 0\cdot\lambda(m_i\sw 1\otimes\ta_i\otimes c_i))
\\
&=& \phi\circ\nu_M(x),
\end{eqnarray*}
where we used that $\phi$ is a right $C$-comodule and right $A$-module
map to derive the second and the third equalities respectively. This
completes the proof
that the functor $-\square_\tC C$ is separable.

Conversely, assume that $-\square_\tC C$ is separable and let $\nu_M$ be
the corresponding splitting of $\Phi_M$.
Define
$$
\lambda : (C\otimes\tA)\square_\tC C\to A, \qquad \lambda =  
(A\otimes\eps_C)\circ\nu_{A\otimes C}(1_A\otimes (C\otimes\tA)\square_\tC
C).
$$
Since $\nu_{A\otimes C}$ is a right $A$-linear map, so is $\lambda$. We
first show that $\nu_M$ can be expressed in terms of $\lambda$. For any 
$M\in \M_A$ and $m\in M$ consider a morphism
$\ell_m: A\otimes C\to M\otimes C$ in $\M_A^C(\psi)$ given by $a\otimes c\mapsto m\cdot
a\otimes c$. Since the splitting of the adjunction $\Phi$ is natural in 
$\M_A^C(\psi)$ we have 
\begin{equation}
\ell_m\circ \nu_{A\otimes C} = \nu_{M\otimes C}\circ
((\ell_m\otimes_A \tA)\square_\tC C).
\label{ell}
\end{equation} 
In particular,
choosing $M=A$ one easily finds that (\ref{ell}) implies that 
$\nu_{A\otimes
C}$ is a left $A$-module map.
Now, if $M\in \M_A^C(\psi)$ one can take the morphism $\rho^M\in {\rm
Hom}_A^C(M, M\otimes C)$, and thus using the naturality of $\nu$, obtain
$
\rho^M\circ\nu_M = \nu_{M\otimes C}\circ (\rho^M\otimes C)$. In view of
(\ref{ell}) this reads for all $\sum_i m_i\otimes\ta_i\otimes c_i\in
(M\otimes \tA)\square_\tC C$ projected down to $(M\otimes_A \tA)\square_\tC C$ 
$$
\rho^M\circ\nu_M(\sum_i m_i\otimes\ta_i\otimes c_i) = \sum_i\ell_{m_i\sw 0}\circ\nu_{A\otimes
C}(1_A\otimes m_i\sw 1\otimes \ta_i \otimes c_i).
$$
Applying $M\otimes \eps_C$ to this last equality and using assumption
that $(f,g)$ is admissible one obtains
$$
\nu_M(\sum_i m_i\otimes \ta_i\otimes c_i) = \sum_i m_i\sw
0\cdot\lambda(m_i \sw 1\otimes\ta_i\otimes c_i).
$$
In particular, the choice $M=A\otimes C$ gives  for all $a\in A$,
$\sum_ic_i\otimes\ta_i\otimes c'_i \in (C\otimes\tA)\square_\tC C$
\begin{equation}
\nu_{A\otimes C}(\sum_i a\otimes c_i\otimes\ta_i \otimes c'_i) = 
a\sum_i\lambda( c_i\sw 2\otimes \ta_i \otimes c'_i)_\alpha\otimes 
c_i\sw 1^\alpha .
\label{nuac}
\end{equation}
We are now ready to show that $\lambda$ satisfies all the conditions of
Definition~\ref{def.int.morphism}. 
Take any $x= \sum_ic_i\otimes\ta_i\otimes c'_i\in (C\otimes \tA)\square_\tC
C$, then
\begin{eqnarray*}
\sum_i \lambda(c_i\otimes\ta_i\otimes c_i'\sw 1) \otimes c_i'\sw 2 
&=& \sum_i (A\otimes \eps_C)\circ
\nu_{A\otimes C} (1_A\otimes c_i\otimes\ta_i \otimes c_i'\sw 1)\otimes
c'_i \sw 2\\
& = & (A\otimes \eps_C\otimes
C)\circ(A\otimes\Delta_C)\circ\nu_{A\otimes C}(1_A\otimes x)\\
& = & \nu_{A\otimes C}(1_A\otimes x)\\
& = & \sum_i\lambda(c_i\sw 1\otimes\ta_i\otimes c_i')_\alpha\otimes c_i
\sw 2^\alpha 
\qquad \mbox{\rm (by (\ref{nuac}))},
\end{eqnarray*}
where we used that $\nu_{A\otimes C}$ is a right $C$-comodule map to
derive the second equality. This proves that $\lambda$ satisfies
(\ref{int.b}). Furthermore, for all $\sum_i c_i\otimes a_i\otimes\ta_i
\otimes c'_i\in (C\otimes A\otimes \tA)\square_\tC C$ we have
\begin{eqnarray*}
\lambda(\sum_ic_i\otimes f(a_i)\ta_i\otimes c'_i) & = &
(A\otimes\eps_C)\circ\nu_{A\otimes C} (\sum_i 1_A\otimes c_i\otimes 
f(a_i)\ta_i\otimes c'_i) \\
&=& (A\otimes\eps_C)\circ\nu_{A\otimes C} (\sum_i (1_A\otimes c_i)\cdot a_i
\otimes \ta_i\otimes c'_i) \\
&=& (A\otimes\eps_C)\circ\nu_{A\otimes C} (\sum_i a_{i\alpha} \otimes
c_i^\alpha \otimes \ta_i\otimes c'_i) \\
&=& \sum_i a_{i\alpha}\lambda(c_i^\alpha \otimes \ta_i\otimes c'_i),
\end{eqnarray*}
where we used the properties of the domain of $\nu_{A\otimes C}$ and
the assumption that $-\square_\tC C$ preserves cokernel of 
$\eq_{M_A\tA}$  to
derive the second equality. This proves that $\lambda$ satisfies
(\ref{int.a}). Finally, for all $c\in C$,
$\Phi_{A\otimes C}(1_A\otimes c) = 1_A\otimes c\sw 1\otimes 1_\tA\otimes
c\sw 2$. Since $\nu_{A\otimes C}$ splits $\Phi_{A\otimes C}$ we have
$1\otimes c = \nu_{A\otimes C}( 1_A\otimes c\sw 1\otimes 1_\tA\otimes
c\sw 2)$. Applying $A\otimes\eps_C$ to this equality one immediately deduces
that  $\lambda$ satisfies (\ref{int.c}).  Therefore the morphism $(f,g)$
is totally integrable. This completes the
proof of the first statement of the theorem.

(2) Given $\z$ as in Definition~\ref{def.coint.morphism}  define for all
$\tM\in \M_A^C(\tpsi)$, $\nu_\tM: \tM\to (\tM\square_\tC C)\otimes_A\tA$,
$\nu_\tM = (\rho_\tM\otimes C\otimes_A
\tA)\circ(\tM\otimes\z)\circ\rho^\tM$. The proof that $\nu_\tM$ is the
required cosplitting is dual to the proof of the corresponding part of
assertion (1). Conversely, given a cosplitting $\nu_\tM$ define $\zeta =
(\eps_\tC\otimes\tA\otimes C\otimes_A \tA)\circ\nu_{\tC\otimes\tA}\circ 
(\tC\otimes 1_\tA)$.
\end{proof}

Notice that the assumption of Theorem~\ref{thm.main}(1) is satisfied if 
${}^\tC C$ is coflat. Dually, the assumption of
Theorem~\ref{thm.main}(2) is satisfied if  ${}_A\tA$ is flat. 
The remainder of the paper is devoted to the analysis of special cases
of Theorem~\ref{thm.main}.

\section{Separable coalgebra-Galois extensions}
\noindent The following notion was introduced 
in \cite{Brz:Fro}. It generalises
the notion of an $H$-integral for a Doi-Hopf datum
\cite[Definition~2.1]{CaeMil:Doi}. 

\begin{definition}
Let $(A,C)_\psi$ be an entwining structure. An {\em integral} in $(A,C)_\psi$ is an
element $\z=\sum_i a_i\otimes c_i\in A\otimes C$ such that for all
$a\in A$, $a\cdot \z = \z\cdot a$. 
Explicitly, we require $\sum_i
aa_i\otimes c_i = \sum_ia_i\psi(c_i\otimes a)$. An integral 
$\z=\sum_i a_i\otimes c_i$ is
said to be {\em normalised} if $\sum_i a_i\eps_C(c_i)=1$.
\label{def.int}
\end{definition}

\begin{example}
Let $A$ be a Hopf algebra and $B\subset A$ be a left $A$-comodule
subalgebra, i.e., a subalgebra of $A$ such that $\Delta_A(B) \subset
A\otimes B$. Consider the coalgebra $C/B^+A$. $C$ is a right $A$-module
in the natural way and there is an entwining structure $(A,C)_\psi$ with
$\psi: c\otimes a\mapsto a\sw 1\otimes c\cdot a\sw 2$. Let $\Lambda \in
C$ be such that for all $a\in A$, $\Lambda\cdot a = \eps_A(a)\Lambda$
and $\eps_C(\Lambda) = 1$. Then $\z = 1\otimes \Lambda$ is an integral
in $(A,C)_\psi$.
\label{ex.homog.int}
\end{example}
\begin{proof}
Clearly, $1_A\eps(\Lambda) = 1_A$. Take any $a\in A$, then
$
(1\otimes\Lambda)\cdot a = a\sw 1\otimes \Lambda\cdot a\sw 2 =
a\otimes\Lambda = a\cdot(1\otimes \Lambda).
$
\end{proof}

In \cite{Brz:Fro} it has been shown that the existence of an integral in
$(A,C)_\psi$ is closely related to the fact that the functor $-\otimes
C:\M_A\to \M_A^C(\psi)$ of Example~\ref{ex.psi.modules}(1) is both left and right adjoint of
the forgetful functor $\M_A^C(\psi)\to \M_A$.  
The following theorem, which is an entwining structure version of
\cite[Theorem~2.14]{CaeMil:sep}, shows that  integrals are closely
related to the separability of $-\otimes C$. 
\begin{theorem}
Let $(A,C)_\psi$ be an entwining structure. The functor $-\otimes C
:\M_A\to \M_A^C(\psi)$ is separable if and only if there exists a
normalised integral in $(A,C)_\psi$.
\label{thm.int}
\end{theorem}
\begin{proof}
Consider an admissible morphism 
$(A,\eps_C) : (A,C)_\psi\to
(A,k)_\twist$. Then $-\otimes C = -\square_kC:\M_A\to \M_A^C(\psi)$. Since ${}_AA$ is
flat, Theorem~\ref{thm.main}(2) can be applied and thus
$-\otimes C$ is separable if and only if $(A,\eps_C)$ is totally
cointegrable, i.e. there exists $\z\in {\rm
Hom}(k, (A\otimes C)\otimes_AA) \cong A\otimes C$ such that conditions
(\ref{coint.a})--(\ref{coint.c}) are satisfied. In this case condition
(\ref{coint.a}) is empty, while condition (\ref{coint.b}) means that
$\z$ is an integral in $(A,C)_\psi$. Finally, condition (\ref{coint.c})
states that $\z$ is normalised.
\end{proof}

The existence of normalised integrals in the canonical entwining structure
associated to a coalgebra-Galois extensions turns out to be equivalent
to the separability of such an extension. First, recall from
\cite{HirSug:sem} 
\begin{definition}
An extension of algebras $B\hookrightarrow A$ is {\em separable} if there exists
$u\in A\otimes_B A$ such that for all $a\in A$,
$au=ua$ and $\mu_{A,B}(u)=1_A$, where $\mu_{A,B}:A\otimes_B A\to A$ is the
natural map induced by $\mu_A$. The element $u$ is called a
{\em separability idempotent}. 
\label{def.sep}
\end{definition}
\begin{proposition}
A coalgebra-Galois extension $A(B)^C$ is separable if and only if there
exists a normalised integral in the canonical entwining structure.
\label{prop.galois.sep}
\end{proposition}
\begin{proof} 
We first show that $\can^{-1}: A\otimes C\to A\otimes_B A$ is an
$(A,A)$-bimodule map, where the $(A,A)$-bimodule structure on $A\otimes
C$ is as in Definition~\ref{def.int}. By construction, $\can^{-1}$ is a left
$A$-module map. For all  
$\z=\sum_i a_i\otimes c_i\in A\otimes C$,  $a\in A$
\begin{eqnarray*}
\can^{-1}(\z\cdot a) &=& \can^{-1}(\sum_i a_ia\sb\alpha\otimes
c_i\sp\alpha) = \sum_i a_i\can^{-1}(a\sb\alpha\otimes
c_i\sp\alpha)\\
& = & \sum_ia_i\can^{-1}(\can(\can^{-1}(1_A\otimes c_i)a))
\qquad\qquad 
\mbox{\rm (def.\ of canonical $\psi$)}\\
& = & \sum_ia_i\can^{-1}(1_A\otimes c_i)a = \can^{-1}(\z)a.
\end{eqnarray*}
Therefore $\z$ is an integral in $(A,C)_\psi$ if and only if for all
$a\in A$, $au=ua$, where $u=\can^{-1}(\z)$. Furthermore, directly from the
definition of the canonical map $\can$, one finds that $(A\otimes \eps_C)\circ \can =
\mu_{A,B}$. Therefore $\z$ is normalised if and only if 
$\mu_{A,B}(u) = 1_A$.
\end{proof}

\begin{example}
In the setting of Example~\ref{ex.homog.int}, view $A$ as a right
$C$-comodule via $\rho^A = (A\otimes\pi)\circ\Delta_A$, where $\pi: A\to
C = A/B^+A$ is the canonical surjection, and assume that $B = \{b\in B
\; | \; \forall a\in A, \; \rho^A(ba) = b\rho^A(a)\}$ (for example, this
holds if either ${}_BA$ or $A_B$ is faithfully flat). Then $A(B)^C$ is a
coalgebra-Galois extension, and if there is $\Lambda \in
C$ such that for all $a\in A$, $\Lambda\cdot a = \eps_A(a)\Lambda$
and $\eps_C(\Lambda) = 1$, then $B\hookrightarrow A$ is separable.
\label{ex.homog.sep}
\end{example}

The introduction of separable extensions in \cite{HirSug:sem} was
motivated by the Hochschild relative homological algebra
\cite{Hoc:rel}. In the case of a coalgebra-Galois extension the
relationship between cohomology and separable extensions can be
expressed in terms of integrals in the canonical entwining structure.
Recall from \cite{Hoc:rel} that if $B$ is a subalgebra of $A$ then for
every $(A,A)$-bimodule $M$ the {\em relative Hochschild cohomology groups}
$H^n(A,B, M)$ are defined as cohomology groups of the complex $(\bigoplus
_{n=0}^\infty C^n(A,B,M), \delta)$, where 
$
C^0(A,B,M) = \{m\in M \; |\; \forall b\in B,\; b\cdot m = m \cdot b\}$,
$$
C^n(A,B,M) = {}_B{\rm
Hom}_B(\underbrace{A\otimes_B A\otimes_B \cdots \otimes_B A}_{\mbox{\rm
$n$-times}},M),  \qquad n>0,
$$ 
and the coboundary 
$ \delta :C^{n}(A,B,M)\to C^{n+1}(A,B,M)$ is given by 
\begin{eqnarray*}
\delta(f)(a_1,\ldots 
,a_{n+1}) &=& a_1\cdot f(a_2,\ldots a_{n+1}) + \sum_{i=1}^{n} (-1)^if(a_1,\ldots
,a_ia_{i+1}, \ldots ,a_{n+1})\\
&& +(-1)^{n+1} f(a_1,\ldots , a_n)\cdot
a_{n+1}. 
\end{eqnarray*}
\begin{corollary}
Let $A(B)^C$ be a coalgebra-Galois extension. Then a normalised
integral in the associated canonical entwining structure exists if and
only if for all $(A,A)$-bimodules $M$, $H^1(A,B,M) = 0$.
\label{cor.sep1}
\end{corollary}
\begin{proof}
By an argument similar to \cite[p.~76]{DeMIng:sep}, one shows, that 
the first relative Hochschild cohomology group
is trivial for all $(A,A)$-bimodules if and only if the extension 
$B\hookrightarrow A$
is separable. Then the assertion follows from
Proposition~\ref{prop.galois.sep}. 
\end{proof}
\begin{corollary}
Let $A(B)^C$ be a coalgebra-Galois extension with a normalised integral
in the canonical entwining structure. Then any $(A,A)$-bimodule which is
semisimple as a $(B,B)$-bimodule is semisimple as an $(A,A)$-bimodule.
\label{cor.sep2}
\end{corollary}
\begin{proof} By Corollary~\ref{cor.sep1}, for all $(A,A)$-bimodules
$M$, $H^1(A,B,M)= 0$. Then 
\cite[Theorem~1]{Hoc:rel} implies the assertion.
\end{proof}

Dually one can consider
\begin{definition}
Let $(A,C,\psi)$ be an entwining structure. A $k$-module map $\y:
C\otimes A\to k$, such that for all $a\in A$, $c\in C$, $c\sw
1\y(c\sw 2\otimes a) = \y(c\sw 1\otimes a_\alpha)c\sw 2^\alpha$
is called a {\em cointegral} in $(A,C)_\psi$. A cointegral $\y$ is
said to be {\em normalised} if $\y\circ(C\otimes 1_A) = \eps_C$.
\label{def.coint}
\end{definition}
\begin{example}
Let $C$ be a Hopf algebra and let $A$ be a right $C$-comodule
algebra. Then $(A,C)_\psi$ is
an entwining structure with $\psi: c\otimes a\mapsto a\sw 0\otimes ca\sw
1$. Let $\kappa\in A^*$ be such that $\kappa(1_A) = 1$ and for all $a\in
A$, $1_C\kappa(a) =
\kappa(a\sw 0)a\sw 1$. Then $\y = \eps_C\otimes \kappa$ is a normalised 
cointegral
in $(A,C)_\psi$.
\label{ex.homog.coint}
\end{example}
\begin{theorem}
Let $(A,C)_\psi$ be an entwining structure. The functor $-\otimes A
:\M^C\to \M_A^C(\psi)$ is separable if and only if there exists a
normalised cointegral in $(A,C)_\psi$.
\label{thm.coint}
\end{theorem}
\begin{proof}
Consider the morphism $(1_A, C): (k,C)_\twist\to (A,C)_\psi$ and apply 
Theorem~\ref{thm.main}(1). 
\end{proof}
\begin{definition}
A coextension of coalgebras $C\twoheadrightarrow B$ is said to be a {\em 
separable coextension}
if there exists a $k$-module map $\upsilon : C\square_B C\to k$ such that
$(C\otimes\upsilon)\circ (\Delta_C\otimes C) = (\upsilon\otimes
C)\circ(C\otimes \Delta_C)$ on $C\square_B C$, and $\upsilon\circ\Delta_C = 
\eps_C$.
\label{def.coseparable}
\end{definition}
\begin{proposition}
An algebra-Galois coextension $C(B)_A$ is separable if and only if there
exists a normalised cointegral in the associated canonical entwining
structure. 
\label{prop.coseparable}
\end{proposition}
\begin{example}
Let $C$ be a Hopf algebra and $A\subset C$ a right comodule
subalgebra of $C$, i.e., $\Delta_C(A)\subset A\otimes C$, so that we are
in the setting of Example~\ref{ex.homog.coint}. Consider the coalgebra  
$B = C/CA^+$, and assume that $A= \{a\in C\; | \; \pi(a\sw 1)\otimes
a\sw 2 = \pi(1_C)\otimes a\}$, where $\pi:C\to B$ is the canonical
surjection (this assumption is satisfied if either ${}_AC$ or $C_A$ 
is faithfully flat).
Then $C\twoheadrightarrow B$ is an $A$-Galois coextension and if there exists
$\kappa\in A^*$ such that for all $a\in A$, $\kappa(a\sw 1)a\sw 2 =
\kappa(a)\eps_C$ and $\kappa(1_A) = 1$, then this coextension is
separable.  
\label{ex.homog.sep*}
\end{example}
When $k = \bf C$, a  rich source of separable coalgebra-Galois 
extensions is provided by 
quantum homogeneous spaces of compact quantum groups \cite{Wor:com}. 
In this case we
are in the setting of Example~\ref{ex.homog.sep*}, with $C$  a compact
quantum group and $A$ a right $C$-homogeneous quantum space. In many
cases $C$ is a faithfully flat right or left $A$-module (see
\cite{MulSch:hom} for examples). The map $\kappa$ is the Haar measure on
$C$ restricted to $A$. Perhaps the simplest example of this situation is
when $C$ is the quantum $SU(2)$ group and $A$ is any of 
the quantum 2-spheres of Podle\'s \cite{Pod:sph}.

\section{Split coalgebra-Galois extensions}
\noindent The following definition is a slightly modified version of
\cite[Definition~4.1]{Brz:Fro}; both definitions describe the same
object if a coalgebra $C$ is a
finitely-generated projective $k$-module.
\begin{definition}
Let $(A,C)_\psi$ be an entwining structure. Any $\gamma\in {\rm
Hom}(C\otimes C, A)$ such that the following diagrams
\begin{equation} 
\begin{CD}
C\otimes C @>{C\otimes \Delta_C}>> C\otimes C\otimes C @>{
\gamma\otimes C}>> A \otimes C \\
@VV{\Delta_C\otimes C}V @.       @|\\
C\otimes C\otimes C @>{C\otimes \gamma}>> C\otimes A @>{\psi}>> A\otimes C
\end{CD}
\label{int.i}
\end{equation}
\begin{equation}
\begin{CD}
C\otimes C\otimes A @>{\gamma\otimes A}>> A\otimes A @>{\mu_A}>> A\\
@VV{C\otimes\psi}V  @.  @AA{\mu_A}A\\
C\otimes A\otimes C @>{\psi\otimes C}>> A\otimes C\otimes C
@>{A\otimes\gamma}>> A\otimes A
\end{CD}
\label{int.ii}
\end{equation}
commute is called an {\em integral map} in $(A,C)_\psi$. An integral map
$\gamma$ is said to be {\em normalised}, if for all $c\in C$,
$\gamma(c\sw 1\otimes c\sw 2) = \eps_C(c)1_A$.
\label{def.integral.map}
\end{definition}
The following theorem is an entwining structure version of
\cite[Theorem~2.3]{CaeMil:sep}.
\begin{theorem}
The forgetful functor $\M_A^C(\psi)\to \M_A$ is separable if and only if
there exists a normalised integral map in $(A,C)_\psi$.
\label{thm.separable}
\end{theorem}
\begin{proof}
Consider an admissible morphism
$(A,\eps_C):(A,C)_\psi\to (A,k)_\twist$. Then $-\otimes_AA:\M_A^C(\psi)
\to \M_A$ is the forgetful functor. 
In this case $(C\otimes A)\square_k C =
C\otimes A\otimes C$, and 
 for all $M\in\M_A$, $\eq_{M_AA} = M$ so that
the assumption of Theorem~\ref{thm.main}(1) holds. Therefore the
forgetful functor is separable if and only if $(A,\eps_C)$ is totally 
integrable, i.e., iff there exists $\lambda\in {\rm Hom}_A(C\otimes
A\otimes C,A)$ satisfying all the conditions of
Definition~\ref{def.int.morphism}. Assume that such a $\lambda$ exists
and define $\gamma = \lambda\circ(C\otimes 1_A\otimes C):C\otimes C\to
A$. Then for all $a\in A$, $c,c'\in C$ we have
\begin{eqnarray*}
a_{\alpha\beta}\gamma(c^\beta\otimes c^{\prime\alpha}) & = &
a_{\alpha\beta}\lambda(c^\beta\otimes 1_A\otimes c^{\prime\alpha}) =
\lambda(c\otimes a_\alpha\otimes c^{\prime\alpha}) \qquad\mbox{\rm (by
(\ref{int.a}))}\\ 
& = & \lambda(c\otimes 1_A\otimes c')a = \gamma(c\otimes c')a,
\end{eqnarray*}
where we used that $\lambda$ is a right $A$-module map to derive the
penultimate equality. Hence the diagram (\ref{int.i}) commutes. Also, (\ref{int.b}) implies
that the diagram (\ref{int.ii}) commutes, while the normalisation of $\gamma$
follows immediately from (\ref{int.c}). Thus we conclude that $\gamma$
is a normalised integral map as required.

Conversely, assume that $\gamma$ is a normalised integral map and define
$\lambda :C\otimes A\otimes C\to A$, $c\otimes a\otimes c' \mapsto
a_\alpha\gamma(c^\alpha\otimes c')$. For all $a,a'\in A$, $c,c'\in C$ we
have
$$
a_\alpha\lambda(c^\alpha\otimes a'\otimes c') = a_\alpha
a'_\beta\gamma(c^{\alpha\beta}\otimes c') =
(aa')_\alpha\gamma(c^\alpha\otimes c') = \lambda(c\otimes aa'\otimes
c'),
$$
where (\ref{diag.A}) was used to obtain the third equality. This proves
that the diagram (\ref{int.a}) commutes. Furthermore
\begin{eqnarray*}
\lambda(c\sw 2\otimes a\otimes c')_\alpha\otimes c\sw 1^\alpha & = &
(a_\delta\gamma(c\sw 2^\delta\otimes c'))_\alpha\otimes c\sw 1^\alpha\\
& = &
a_{\delta\alpha}\gamma(c\sw 2^\delta\otimes c')_\beta\otimes c\sw
1^{\alpha\beta} \qquad \qquad \mbox{\rm (by (\ref{diag.A}))}\\
& = & a_{\alpha}\gamma(c^\alpha\sw 2\otimes c')_\beta\otimes c^\alpha \sw
1^{\beta} \qquad\qquad  \mbox{\rm (by (\ref{diag.B}))}\\ 
& = & a_{\alpha}\gamma(c^\alpha\otimes c'\sw 1)_\beta\otimes c'\sw
2 \qquad \qquad\mbox{\rm (by (\ref{int.ii}))}\\
& = & \lambda(c\otimes a\otimes c'\sw 1)\otimes c'\sw 2.
\end{eqnarray*}
This proves that diagram (\ref{int.b}) commutes. Also,
\begin{eqnarray*}
\lambda(c\otimes aa'_\alpha\otimes c^{\prime\alpha}) & = &
(aa'_\alpha)_\beta \gamma (c^\beta\otimes c^{\prime\alpha}) = a_\beta
a'_{\alpha\delta}  \gamma (c^{\beta\delta}\otimes c^{\prime\alpha})
\qquad \mbox{\rm (by (\ref{diag.A}))}\\
& = & a_\beta\gamma(c^\beta\otimes c')a' = \lambda(c\otimes a\otimes a')
\qquad\qquad \quad  \mbox{\rm (by (\ref{int.i}))}.
\end{eqnarray*}
Therefore $\lambda$ is a right $A$-module map, and, consequently the
morphism $(A,\eps_C)$ is integrable. The fact that it is totally
integrable follows immediately from the normalisation of $\gamma$. 
\end{proof} 
\begin{example}
Let $(A,C)_\psi$ be the canonical entwining structure associated to a
pointed algebra-Galois coextension $C(k)_A$ of $k$. Then the forgetful
functor $\M_A^C(\psi)\to \M_A$ is separable.
\label{ex.separable}
\end{example}
\begin{proof}
Since $B=k$, $C\square_B C = C\otimes C$, and  we define 
 $\gamma = (\eps_C\otimes A)\circ \cocan^{-1}:C\otimes C\to A$. We show
that $\gamma$ is a normalised integral map. First notice that since 
$\cocan^{-1}$ is a left $C$-comodule map, one has $\cocan^{-1} =
(C\otimes\gamma)\circ(\Delta_C\otimes C)$. Applying the definition of
the canonical entwining map in Example~\ref{cocan.ex} to $\cocan^{-1}$
one immediately obtains $\psi\circ \cocan^{-1} = 
(\gamma\otimes C)\circ(C\otimes\Delta_C)$, i.e.
$\psi\circ(C\otimes\gamma)\circ(\Delta_C\otimes C) = (\gamma\otimes
C)\circ(C\otimes\Delta_C)$. Thus we conclude that $\gamma$ satisfies
condition (\ref{int.i}).

Let $\kappa:A\to k$ be the algebra map making $C(k)_A$ a pointed
algebra-Galois coextension. One easily finds that $\rho_C =
(\kappa\otimes C)\circ \psi$ and $C\otimes \kappa = (C\otimes\eps_C)\circ
\cocan$. The map $\gamma$ is the {\em cotranslation map}, so, as explained in 
\cite[Theorem~3.5]{BrzHaj:coa}, it has the following properties
\begin{equation}
\mu_A\circ(\gamma\otimes A) = \gamma\circ(C\otimes\rho_C),
\label{cot.2}
\end{equation}
\begin{equation}
\mu_A\circ(\gamma\otimes\gamma)\circ(C\otimes \Delta_C\otimes C) =
\gamma\circ(C\otimes\eps_C\otimes C).
\label{cot.3}
\end{equation}
Using all these properties we obtain
\begin{eqnarray*}
\mu_A\circ(\gamma\otimes A) & = &
\gamma\circ(C\otimes\rho_C)\qquad\qquad\qquad\qquad\qquad
\qquad\qquad\qquad\qquad  \mbox{\rm (by
(\ref{cot.2}))}\\ 
& = & \gamma\circ(C\otimes\kappa\otimes C)\circ(C\otimes\psi) \\
& = & \gamma\circ(C\otimes\eps_C\otimes  C)\circ(\cocan\otimes
C)\circ(C\otimes\psi)\\ 
& = & \mu_A\circ(\gamma\otimes\gamma)\circ(C\otimes \Delta_C\otimes
C)\circ(\cocan\otimes C)\circ (C\otimes\psi) \quad \mbox{\rm (by
(\ref{cot.3}))}\\ 
& = & \mu_A\circ(A\otimes\gamma)\circ(\psi\otimes C)\circ (C\otimes\psi)
\qquad\qquad\qquad \qquad \mbox{\rm (def. of $\psi$)}.
\end{eqnarray*}
This proves that $\gamma$ is an integral map. Finally, $\gamma$ is
normalised by the normalisation property of the cotranslation map (cf.
\cite[Theorem~3.5]{BrzHaj:coa}). 
\end{proof}

As explained in \cite{CaeMil:sep} the separability of the forgetful
functor implies various Maschke-type theorems. Thus, similarly as in
\cite{Brz:Fro} we have
\begin{corollary}
If there is a normalised integral map in $(A,C)_\psi$, then 

(1) Every 
object in $\M_A^C(\psi)$ which is semisimple as an object in $\M_A$ is
semisimple as an object in $\M_A^C(\psi)$. 

(2) Every object in
$\M_A^C(\psi)$ which is projective (resp.\ injective) as a 
right $A$-module is
a projective (resp.\ injective) object in $\M_A^C(\psi)$. 

(3) $M\in \M_A^C(\psi)$ is
projective as a right $A$-module if and only if there exists $V\in \M^C$
such that $M$ is a direct summand of $V\otimes A$ in $\M_A^C(\psi)$
($V\otimes A$ is an entwined module by Example~\ref{ex.psi.modules}(2)).
\end{corollary}
In the case of a coalgebra-Galois extension, 
the existence of normalised
integral maps in the canonical entwining structure is closely related 
to the coalgebra-Galois extension being
a split extension. Recall from
\cite{Pie:ass}\cite{Kad:Jon}
\begin{definition}
An extension of algebras $B\hookrightarrow A$ is called a 
{\em split extension} if there exists a
unital $(B,B)$-bimodule map $E:A\to B$. The map $E$ is called a {\em
conditional expectation}.
\label{def.split}
\end{definition}
\begin{proposition}
A coalgebra-Galois extension $A(B)^C$ is a split extension if and only 
if there
exists $\phi\in{\rm Hom}(C,A)$ such that

(i) $\forall c\in C, \qquad \psi(c\sw 1\otimes \phi(c\sw 2)) =
\phi(c)\rho^A(1_A)$,

(ii)  $\sum_i a^i\phi(c_i) = 1_A$, where $\sum_i a^i\otimes c_i =
\rho^A(1_A)$.

(iii) $\forall b\in B,\; c\in C,\qquad  b_\alpha\phi(c^\alpha) = 
\phi(c)b$.
\label{proposition.split}
\end{proposition}
\begin{proof}
As explained in the proof of \cite[Proposition~4.4]{Brz:mod}, given a
unital $(B,B)$-bimodule map $E:A\to B$ there exists $\phi\in{\rm
Hom}(C,A)$ satisfying conditions (i)--(iii). Explicitly, $\phi =
(A\otimes_B E)\circ \can^{-1}\circ (1_A\otimes C)$. Conversely, given
$\phi\in{\rm Hom}(C,A)$ satisfying (i), \cite[Theorem~4.3]{Brz:mod}
implies that $E:A\to B$, $a\mapsto a\sw 0\phi(a\sw 1)$ is a left
$B$-module map. Clearly, condition (ii) implies $E$ is unital.
Furthermore, for all $a\in A$, $b\in B$
$$
E(ab) = (ab)\sw 0\phi((ab)\sw 1) = a\sw 0b_\alpha\phi(a\sw 1^\alpha) =
a\sw 0\phi(a\sw 1) b = E(a)b,
$$
where we used that $A\in \M_A^C(\psi)$ and the assumption (iii) to derive
the second and third equalities respectively. This completes the proof.
\end{proof}

Let $(A,C)_\psi$ be an entwining structure and assume that
$A\in\M_A^C(\psi)$. Define $B$ as in Example~\ref{can.ex}. Then 
one can consider a covariant functor $(-)_0: \M_A^C(\psi)
\to \M_{B}$
$$
M\mapsto M_0 := \{m\in M\; | \; \forall a\in A, \;\; \rho^M(m\cdot a) =
m\rho^A(a)\}. 
$$
Notice, in particular, that $B=A_0$. 
As explained in \cite{Brz:mod} the functor $(-)_0$ is the right adjoint
of the functor $-\otimes_BA:\M_B\to \M_A^C(\psi)$.
\begin{corollary}
If a coalgebra-Galois extension $A(B)^C$ is a split extension then 
${}_BA$ is a
faithfully flat module. Consequently, the functors $-\otimes_B A:\M_B\to
\M_A^C(\psi)$ and $(-)_0 :\M_A^C(\psi)\to \M_B$ are inverse
equivalences.
\label{cor.split}
\end{corollary}
\begin{proof}
The first assertion follows from \cite[Proposition~4.4]{Brz:mod}, while
the second is the consequence of \cite[Corollary~3.11]{Brz:mod}. 
\end{proof}

\begin{proposition}
Let $(A,C)_\psi$ be the canonical entwining structure associated to a
coalgebra-Galois extension $A(B)^C$. If there is a normalised integral
map in $(A,C)_\psi$ then $B\hookrightarrow A$ is a split extension.
\label{prop.split}
\end{proposition}
\begin{proof}
Let $\gamma: C\otimes C\to A$ be a normalised integral map in $(A,C)_\psi$,
and take $\phi:C \to A$, $c\mapsto \sum_i
a^i_{\alpha}\gamma(c^\alpha\otimes c_i)$, where
$\sum_ia^i\otimes c_i = \rho^A(1_A)$. Notice that the fact that $A$ is
an $(A,C)_\psi$ module implies that for all $a\in A$, $\rho^A(a) =
\sum_ia^ia_\alpha\otimes c_i^\alpha$. Furthermore, since $\rho^A$ is a
coaction  we have
\begin{equation}
\sum_{i,j}a^ja^i_\alpha\otimes c_j^\alpha \otimes c_i = \sum_i a^i\otimes c_i\sw
1\otimes c_i\sw 2.
\label{star}
\end{equation}
We now show that $\phi$ satisfies all the conditions of
Proposition~\ref{proposition.split}. For all $c\in C$
\begin{eqnarray*}
\psi(c\sw 1\otimes \phi(c\sw 2)) & = &  (\sum_i
a^i_{\alpha}\gamma(c\sw 2^\alpha\otimes c_i))_\beta\otimes c\sw
1^\beta \\
& = &\sum_ia^i_{\alpha\delta}\gamma(c\sw 2^\alpha\otimes c_i)_\beta\otimes c\sw
1^{\delta\beta} \qquad\qquad \mbox{\rm (by (\ref{diag.A}))}\\
& = & \sum_ia^i_{\alpha}\gamma(c^\alpha \sw 2\otimes c_i)_\beta\otimes 
c^{\alpha}\sw 1^\beta  \qquad\qquad \mbox{\rm (by (\ref{diag.B}))}\\ 
& = & \sum_ia^i_{\alpha}\gamma(c^\alpha \otimes c_i\sw 1)\otimes 
c_i\sw 2 \qquad\qquad \mbox{\rm (by (\ref{int.i}))}\\ 
& = &    \sum_{i,j}(a^ja^i_{\beta})_{\alpha}\gamma(c^\alpha \otimes 
c_j^\beta)\otimes 
c_i \qquad\qquad \mbox{\rm (by (\ref{star}))}\\ 
& = &    \sum_{i,j}a^j_\alpha a^i_{\beta\delta}\gamma(c^{\alpha\delta}
\otimes 
c_j^\beta)\otimes 
c_i \qquad\qquad \mbox{\rm (by (\ref{diag.A}))}\\ 
& = &  \sum_{i,j}a^j_\alpha\gamma(c^\alpha\otimes c_j)a^i\otimes c_i 
= \phi(c)\rho^A(1_A) 
\qquad \mbox{\rm (by (\ref{int.ii}))}  
\end{eqnarray*}
Using normalisation of $\gamma$ as well as (\ref{star}) one easily finds
that $\sum_i a^i\phi(c_i) = 1_A$. Finally, take any $b\in B$, $c\in C$ and
compute
\begin{eqnarray*}
b_\alpha\phi(c^\alpha) & = & \sum_i b_\alpha 
a^i_{\beta}\gamma(c^{\alpha\beta}\otimes c_i)
=  \sum_i (b 
a^i)_\alpha\gamma(c^\alpha\otimes c_i) \qquad \mbox{\rm
(by (\ref{diag.A}))} \\
& = &  \sum_i ( 
a^ib_\beta )_\alpha\gamma(c^\alpha\otimes c_i^\beta ) =   \sum_i  
a^i_\alpha b_{\beta\delta}\gamma(c^{\alpha\delta}\otimes c_i^\beta )\qquad \mbox{\rm
($b\in B$,  (\ref{diag.A}))} \\
& = &   \sum_i  
a^i_\alpha \gamma(c^{\alpha}\otimes c_i )b = \phi(c)b \qquad \qquad 
\qquad\qquad \mbox{\rm
(by (\ref{int.ii}))} 
\end{eqnarray*} 
Therefore $\phi$ satisfies all the conditions of Proposition~\ref{proposition.split}
and, consequently, $B\hookrightarrow A$ is a split extension.
\end{proof}

Dually to Definition~\ref{def.integral.map} we can consider
\begin{definition}
Let $(A,C)_\psi$ be an entwining structure. Any $\zeta\in {\rm
Hom}(C, A\otimes A)$ such that the following diagrams
\begin{equation}
\begin{CD}
C\otimes A@>{\zeta\otimes A}>> A\otimes A \otimes A@>{A\otimes
\mu_A}>> A\otimes A \\
@| @.       @AA{\mu_A\otimes A}A\\
C\otimes A @>{\psi}>> A\otimes C @>{A\otimes\zeta}>> A\otimes A\otimes A
\end{CD}
\label{coint.i}
\end{equation}
\begin{equation}
\begin{CD}
C @>{\Delta_C}>> C\otimes C @>{\zeta\otimes C}>> A\otimes A\otimes C\\
@VV{\Delta_C}V  @.  @AA{A\otimes\psi}A\\
C\otimes C @>{C\otimes\zeta}>> C\otimes A\otimes A
@>{\psi\otimes A}>> A\otimes C\otimes A
\end{CD}
\label{coint.ii}
\end{equation}
commute is called a {\em cointegral map} in $(A,C)_\psi$. A cointegral map
$\zeta$ is said to be {\em normalised}, if 
$\mu_A\circ \zeta= 1_A\circ \eps_C$.
\label{def.cointegral.map}
\end{definition}
\begin{theorem}
The forgetful functor $\M_A^C(\psi)\to \M^C$ is separable if and only if
there exists a normalised cointegral map in $(A,C)_\psi$.
\label{thm.separable*}
\end{theorem}
\begin{proof}
Consider an admissible morphism 
$(1_A,C): (k,C)_\sigma\to
(A,C)_\psi$ and apply Theorem~\ref{thm.main}(2).
\end{proof}
\begin{example}
Let $(A,C)_\psi$ be a canonical entwining structure associated to a
copointed coalgebra-Galois extension $A(k)^B$ of $k$. Then the forgetful
functor $\M_A^C(\psi)\to \M^C$ is separable. 
\label{ex.separable*}
\end{example}
In this case a normalised cointegral map 
is $\zeta = \can^{-1}\circ (1_A\otimes C)$.

\section{Strongly separable coalgebra-Galois extensions}
\noindent In this  section we combine the results of previous two sections to
determine when a coalgebra-Galois extension is a strongly separable
extension. Such  an extension was introduced in \cite{Kad:Jon} in order
to describe algebraic aspects of the
Jones knot polynomial.
\begin{definition}
An extension of algebras $B\hookrightarrow A$ is called a {\em strongly separable
extension} if it is a separable and split extension, and there exist a
separation idempotent $u = \sum_iu_i\otimes u^i$, a conditional
expectation $E:A\to B$ and a unit $\tau\in k$ such that for all $a\in
A$,

(i) $\sum_i E(au_i)u^i = a\tau$

(ii) $\sum_i u_iE(u^ia) = a\tau$.
\label{def.str}
\end{definition}
\begin{proposition}
Let $A(B)^C$ be a coalgebra-Galois extension. If there exist 
a normalised integral $\z=\sum_ia_i\otimes c_i$ and
a normalised integral map $\gamma\in {\rm Hom}(C\otimes C,A)$ in the
canonical entwining structure $(A,C)_\psi$, and a unit $\tau\in k$ such
that\\
\indent (i) $
\sum_i a_i1_A\sw 0_\alpha\gamma(c_i^\alpha\otimes 1_A\sw 1) = \tau,
$\\
\indent (ii) $\sum_ia_i\sw 0\gamma(a_i\sw 1\otimes c_i) = \tau ,$\\
then $B\hookrightarrow A$ is a strongly separable extension.
\label{prop.str}
\end{proposition}
\begin{proof}
By Proposition~\ref{prop.galois.sep}, $B\hookrightarrow A$ is separable 
with
$u=\sum_i u_i\otimes u^i = \can^{-1}(\z)$, while by
Proposition~\ref{prop.split}, $B\hookrightarrow A$ is split with a
conditional expectation
$E:A\mapsto a\sw 01_A\sw 0_\alpha\gamma(a\sw 1^\alpha\otimes 1_A\sw 1) =
(a1_A\sw 0)\sw 0\gamma((a1_A\sw 0)\sw 1\otimes 1_A\sw 1)$. Take any $a\in A$
and compute:
\begin{eqnarray*}
\sum_i E(au_i)u^i & = & \sum_i E(u_i)u^i a \qquad \qquad \qquad \mbox{\rm ($u$
is an integral)}\\
& = & \sum_i u_i\sw 0 1_A\sw 0_\alpha \gamma(u_i\sw 1^\alpha\otimes
1_A\sw 1)u^ia\\
& = & \sum_i u_i\sw 0 1_A\sw 0_\alpha u^i_{\beta\delta}
\gamma(u_i\sw 1^{\alpha\delta}\otimes
1_A\sw 1^\beta)a \qquad \mbox{\rm (by (\ref{int.ii}))}\\
& = & \sum_i u_i\sw 0 (1_A\sw 0  u^i_{\beta})_\alpha
\gamma(u_i\sw 1^{\alpha}\otimes
1_A\sw 1^\beta)a \qquad \mbox{\rm (by (\ref{diag.A}))}\\
& = & \sum_i u_i\sw 0 u^i\sw 0_\alpha
\gamma(u_i\sw 1^{\alpha}\otimes
u^i\sw 1)a \qquad \mbox{\rm ($A\in \M_A^C(\psi)$)}\\
& = & \sum_i (u_i u^i\sw 0)\sw 0
\gamma((u_i u^i\sw 0)\sw 1\otimes
u^i\sw 1)a \qquad \mbox{\rm ($A\in \M_A^C(\psi)$)}\\
& = & \sum_i a_i\sw 0\gamma(a_i\sw 1\otimes c_i)a =\tau a \qquad \mbox{ \rm
($\z=\can(u)$)} 
\end{eqnarray*}
Therefore the condition Definition~\ref{def.str}(i) is satisfied. 
Furthermore
\begin{eqnarray*}
\sum_i u_iE(u^ia) & = & \sum_i au_iE(u^i) \qquad \qquad
\qquad \mbox{\rm ($u$
is an integral)}\\
& = & \sum_i au_i u^i\sw 01_A\sw 0 _\alpha\gamma(u^i\sw 1^\alpha
\otimes 1_A\sw 1)\\
& = & \sum_i aa_i1_A\sw 0_\alpha\gamma(c_i^\alpha\otimes 1_A\sw 1) 
= \tau a \qquad \mbox{ \rm ($\z=\can(u)$)}
\end{eqnarray*}
This proves Definition~\ref{def.str}(ii) and thus  completes the proof of
the proposition. 
\end{proof}
\begin{proposition}
Let $k$ be a field and let $A(B)^C$ be a coalgebra-Galois extension with
both $A$ and $B$ finite dimensional. Suppose that $A_B$ is free. Then $B\hookrightarrow A$ is a
strongly separable
extension if and only if there exists a normalised  integral $\z = \sum_i a_i\otimes
c_i$ in the canonical entwining structure $(A,C)_\psi$, a map $\phi:C\to
A$ satisfying conditions (i)--(iii) in Proposition~\ref{proposition.split}, and a
non-zero $\tau\in k$ such that 
\begin{equation}
\sum_ia_i\phi(c_i) = \tau.
\label{f.dim}
\end{equation}
\label{prop.f.dim}
\end{proposition}
\begin{proof}
By Proposition~\ref{prop.galois.sep}, $B\hookrightarrow A$ is 
separable with
$u=\sum_i u_i\otimes u^i = \can^{-1}(\z)$, 
while by
Proposition~\ref{proposition.split}, $B\hookrightarrow A$ is split with a conditional
 expectation
$E:a\mapsto a\sw 0\phi(a\sw 1)$. By \cite[Remark~1.4(d)]{FisMon:Fro},
Definition~\ref{def.str}(i) holds provided that condition
Definition~\ref{def.str}(ii) holds. Thus it
suffices to prove that (\ref{f.dim}) is a sufficient and necessary
condition for Definition~\ref{def.str}(ii). Take any $a\in A$
and compute:
\begin{eqnarray*}
\sum_i u_iE(u^ia) & = & \sum_i au_iE(u^i) \qquad \qquad
\qquad \mbox{\rm ($u$
is an integral)}\\
& = & \sum_i au_iu^i\sw 0\phi(u^i\sw 1) = \sum_i
aa_i\phi(c_i).
\end{eqnarray*}
Therefore $\sum_i u_iE(u^ia) = \tau a$ if and only if
(\ref{f.dim}) holds. 
\end{proof}

\end{document}